\documentclass[12pt]{article}
\usepackage[english]{babel}
\usepackage{amsmath}
\usepackage{amssymb}
\usepackage{amsthm}

        \newtheorem{lemma}{Lemma}[section]
        
        \newtheorem{theorem}[lemma]{Theorem}
        \newtheorem{corollary}[lemma]{Corollary}

        \newtheorem{remark}[lemma]{Remark}

\newcommand{\ov}{\overline}
\newcommand{\bear}{\begin{array}}
\newcommand{\enar}{\end{array}}
\newcommand{\beq}{\begin{equation}}
\newcommand{\eeq}{\end{equation}}
\newcommand{\beqn}{\begin{eqnarray}}
\newcommand{\eeqn}{\end{eqnarray}}
\newcommand{\beit}{\begin{itemize}}
\newcommand{\eeit}{\end{itemize}}

\newcommand{\rsp}{{\bf R}}

\newcommand{\ds}{\displaystyle}

\newcommand{\g}{\gamma}
\renewcommand{\l}{\lambda}

\newcommand{\s}{\sigma}

\newcommand{\Om}{\Omega}
\renewcommand{\phi}{\varphi}
\renewcommand{\a}{\alpha}
\renewcommand{\b}{\beta}

\newcommand{\wtil}[1]{\widetilde{#1}}
\newcommand{\pn}{\par \noindent}
\newcommand{\med}{\medskip}
\newcommand{\qq}{\qquad}
\newcommand{\q}{\quad}
\newcommand{\hookto}{\hookrightarrow}

\newcommand{\media}[1]{\kern 0.4ex {-} \kern -2.0 ex {\int}_{\kern -0.9 ex
 {#1}}\;}

\newcommand{\PI}{\textrm{P}(\textrm{K}), \textrm{K}\in\{\textrm{D,N}\}}
\newcommand{\PII}{\textrm{P}(\textrm{K})}
\newcommand{\EPI}{\emph{P}(\emph{K})$, $\emph{K}\in\{\emph{D,N}\},}

\title{Parabolic integrodifferential identification \\
problems related to radial memory kernels II\footnote{Work partially
supported by the Italian Ministero dell'Universit\`a e della Ricerca
Scientifica e Tecnologica (M.U.R.S.T.).}}

\author{ Alberto Favaron (Milan), Alfredo Lorenzi (Milan)
\footnote{The authors are members of G.N.A.M.P.A. of the Italian
Istituto Nazionale di Alta Matematica (I.N.d.A.M.)}}

\date{}

\begin{document}
\maketitle
\pn
{\bf Abstract.} We are concerned with the problem of recovering the radial kernel $k$,
depending also on time, in the parabolic integro-differential equation
$$D_{t}u(t,x)={\cal A}u(t,x)+\int_0^t\!\! k(t-s,|x|)\mathcal{B}u(s,x)ds
+\int_0^t\!\! D_{|x|}k(t-s,|x|)\mathcal{C}u(s,x)ds+f(t,x),$$
${\cal A}$ being a uniformly elliptic second-order linear
operator in divergence form. We single out a special class of operators ${\cal A}$
and two pieces of suitable additional information for which the
problem of identifying $k$ can be uniquely solved locally in time
when the domain under consideration is a ball or a disk.
\med \pn
\med \pn

\section{Posing the identification problem}
\setcounter{equation}{0}
The present paper is strictly related to our previous one \cite{3}.
Indeed, the
problem we are going to investigate consists, as in \cite{3},
 in identifying an unknown radial
memory kernel $k$ also depending on time, which appears in the following
integro-differential equation related to the ball
$\Omega\!=\!\{x\!\!=\!\! (x_1,x_2,x_3)\in\mathbb{R}^3\!:\!|x|<R\}$,
 $R>0$ and $|x|={(x_1^2+ x_2^2+ x_3^2)}^{\!1/2}$:
\begin{eqnarray}\label{problem}
D_{t}u(t,x)=\mathcal{A}u(t,x)+\!\int_0^t\!\! k(t-s,|x|)\mathcal{B}u(s,x)ds+\!
\int_0^t\!\! D_{|x|}k(t-s,|x|)\mathcal{C}u(s,x)ds\;+\!\!\!\!&f(t,x),&
\nonumber\\[2mm]\hskip 8truecm \forall\, (t,x)\in [0,T]
\times\Omega. & &
\end{eqnarray}
We emphasize that the aim of the present paper is to study the identification
problem related to $(\ref{problem})$ when the domain $\Om$ is a {\emph{full}}
ball. This is exactly a singular domain for our problem as we noted in Remark 2.9 in
\cite{3}, where we were able to recover the kernel $k$ only
in the case of a spherical
corona or an annulus $\Om$. In this paper we show that our identification
problem can actually be solved in suitable weighted spaces if we appropriately
restrict the class of admissible differential operators $\cal{A}$ to a class
whose coefficients have an appropriate structure in a neighbourhood of the
centre $x=0$ of $\Om$, which turns out to be a ``singular point'' for our
problem.\\
In equation $(\ref{problem})\;\mathcal{A}$ and $\mathcal{B}$ are two
second-order linear differential operators, while $\mathcal{C}$ is a
first-order differential operator having  the following forms,
respectively:
\begin{eqnarray}
\label{A}
\mathcal{A}=\sum_{j=1}^{3}D_{x_j}\big(\sum_{k=1}^{3}a_{j,k}(x)D_{x_k}
\big),\ \
\label{B}\mathcal{B}=\sum_{j=1}^{3}D_{x_j}\big(\sum_{k=1}^{3}b_{j,k}(x)
D_{x_k}
\big),\ \
\label{C}\mathcal{C}=\sum_{j=1}^{3}c_{j}(x)D_{x_j}.
\end{eqnarray}
In addition, operator $\mathcal{A}$ has a very special structure,
since its coefficients $a_{i,j}$, $i,j=1,2,3,$ have the following particular
representation,
(cf. \cite{3}, formula $(2.4)$, where $(b,d)$ is changed
in $(-b, -d)$):
\begin{equation}\label{condsuaij}
\left\{\!\!\!\begin{array}{lll}
 a_{1,1}(x)\!\!\!&=&\!\!\!a(|x|)+
\displaystyle\frac{(x_2^2+x_3^2)[c(x)+b(|x|)]}{|x|^2}-
\frac{x_1^2d(|x|)}{|x|^2},\\[5mm]
 a_{2,2}(x)\!\!\!&=&\!\!\!a(|x|)+
\displaystyle\frac{(x_1^2+x_3^2)[c(x)+b(|x|)]}{|x|^2}-
\frac{x_2^2d(|x|)}{|x|^2},
\\[5mm]
 a_{3,3}(x)\!\!\!&=&\!\!\!a(|x|)+
\displaystyle\frac{(x_1^2+x_2^2)[c(x)+b(|x|)]}{|x|^2}
-\frac{x_3^2d(|x|)}{|x|^2},
\\[5mm]
 a_{j,k}(x)\!\!\!&= &\!\!\!a_{k,j}(x)=\displaystyle
-\frac{x_jx_k[b(|x|)+c(x)+d(|x|)]}{|x|^2},\qq
1\le j,k\le 3,\ j\neq k,
\end{array}\right.
\end{equation}
\pn
where the functions $a$, $b$, $c$, $d$  are {\it non-negative} and
enjoy the following properties:
\begin{eqnarray}\label{abcd}
\label{regular}&a, b, d \in C^{2}\big([0,R]\big)\,,\q\;
c\in C^2(\ov\Om),\qq&\\[1,7mm]
\label{bcd} &a(r)> d(r)\,,\q \forall r\in [0,R]\,,
\q b(0)+c(0)=0\,,\q d(0)=0 .&
\end{eqnarray}
In particular, we note that each coefficient $a_{i,j}$ is Lipschitz-continuous
in ${\ov \Om}$.
\pn
We now introduce the function $h$ defined by
\begin{equation}\label{H}
h(r)=a(r)-d(r),\qq\forall\,r\in [0,R]\,,
\end{equation}
and which is non-negative by virtue of $(\ref{bcd})$. Then,
as we noted in \cite{3}, for every $x\in\ov\Omega$
and $\xi\in\mathbb{R}^3$ we have
\begin{eqnarray}\label{unel1}
 \sum_{j,k=1}^{3}a_{j,k}(x){\xi}_j{\xi}_k
\!\!&\geqslant&\!\! a(|x|){|\xi|}^2+\frac{b(|x|)+c(x)}{|x|^2}\,
{|x\wedge\xi |}^2-\frac{d(|x|)}{|x|^2}\,{[x\cdot\xi]}^2\nonumber\\[2mm]
&\geqslant&\!\! a(|x|){|\xi|}^2+\frac{b(|x|)}{|x|^2}\,
{|x\wedge\xi |}^2-\frac{d(|x|)}{|x|^2}\,{[x\cdot\xi]}^2\geqslant
 h(|x|)\,|\xi|^2\ge 0,\qq\q
\end{eqnarray}
where $ \wedge$ and $\cdot$ denote, respectively, the wedge and inner
products in $\mathbb{R}^3$.\\
>From $(\ref{unel1})$ it follows that the condition
of uniform ellipticity of $\cal{A}$, i.e.
\begin{equation}\label{unel}
{\alpha}_1|\xi{|}^2\leqslant\sum_{j,k=1}^{3}a_{j,k}(x){\xi}_j{\xi}_k
\leqslant{\alpha}_2|\xi{|}^2,\qquad\,\forall\,
(x,\xi)\in\Omega\times \mathbb{R}^3\;,
\end{equation}
 is trivially satisfied with
$\a_1\!=\!\min_{r\in [0,R]}h(r)$ and
$\a_2\!=\!\|h+b\|_{C([0,R])}+\|c\|_{C(\ov\Om)}$.\\
Then we prescribe the {\it{initial condition}}:
\begin{equation} \label{u0}
u(0,x)=u_0(x)\,,\;\;\;\;\qquad \forall\, x\in\Omega\,,
\end{equation}
 $u_0:\overline{\Omega}\rightarrow\mathbb{R}$ being a given smooth function,
as well as one of the following boundary value conditions,
where $u_1\!:\![0,T]\!\times\!\overline{\Omega}\!\rightarrow\!\mathbb{R}$
is a given smooth function:
 \begin{alignat}{2}
\label{D11}
& (\textrm{D})\quad\qquad & u(t,x)=u_1(t,x),\qquad\qquad \quad
& \forall\, (t,x)\in [0,T]\times\partial\mbox{}\Omega,\,\\[2mm]
\label{N11}
& (\textrm{N})\quad\qquad & \frac{\partial u}{\partial n}(t,x)
= \frac{\partial u_1}{\partial n}(t,x),\qquad\quad \quad
& \forall\, (t,x)\in [0,T]\times\partial\mbox{}\Omega.
\end{alignat}
Here D and N stand, respectively, for  the Dirichlet and Neumann
boundary conditions, whereas  $n$ denotes the outwarding normal to
$\partial\mbox{}\Omega$.
\begin{remark}\label{conormal}
\emph{The conormal vector associated with the matrix $\{a_{j,k}(x)\}_
{j,k=1}^{3}$ defined by $(\ref{condsuaij})$ and the boundary $\partial\mbox{}
\Om$ coincides with $R^{-1}[a(R)-d(R)]x$, i.e. with the outwarding
normal $n(x)$.}
\end{remark}
To determine the radial memory kernel $k$ we need also the two following
pieces of information:
\begin{eqnarray}
\label{g11}\!\!\!&\Phi&\!\!\!\!\![u(t,\cdot)](r)\!:= g_1(t,r),\,\qquad
\forall\,(t,r)\in[0,T]\times (0,R),\\[1,5mm]
\label{g22}\;\!\!\!&\Psi&\!\!\!\!\![u(t,\cdot)]\!:= g_2(t),\;\,\q\quad
\q\q\forall\,t\in[0,T],
\end{eqnarray}
where, representing with $(r, \varphi, \theta)$ the usual spherical
co-ordinates with pole at $x=0$,  $\Phi$ and $\Psi$ are two linear operators
 acting, respectively, on the angular variables
$\varphi,\,\theta$ only and all the space variables $r,\,\varphi,\,\theta$.\\
\vskip -0,3truecm
\pn{\it{Convention:}} from now on we will denote by
$\textrm{P}(\textrm{K}),\,\textrm{K}\in\{\textrm{D,N}\}$, the
identification problem consisting of $(\ref{problem}), (\ref{u0})$,
the boundary condition $(\textrm{K})$ and $(\ref{g11}),(\ref{g22})$.
\\
\vskip -0,3truecm

\pn An example of admissible linear operators $\Phi$ and $\Psi$ is the
following:
\begin{align}\label{Phi1}
&\Phi [v](r):= \int_{\!0}^{\pi}\!\!\sin\!\theta\mbox{}d\theta\!
\int_{\!0}^{2\pi}\!\!\!\!v(rx') d\varphi\;, &\\[3mm]
\label{Psi1}
&\Psi[v]:=\int_{\!0}^{R}\!\!r^2 dr\!\!\int_{\!0}^{\pi}\!\!
\sin\!\theta\mbox{} d\theta\!
\int_{\!0}^{2\pi}\!\!\!\!\psi(rx')v(rx')  d\varphi \;\,,&
\end{align}
where $x'\!=\!(\cos\!\varphi\sin\!\theta,\,\sin\!\varphi\sin\!\theta,\,
\cos\!\theta) $, while $\psi:\overline{\Omega}\rightarrow\mathbb{R}$
is a smooth assigned function.
\begin{remark}
\emph{We note that $(\ref{Phi1})$ coincides with $(1.12)$ in \cite{3}
with $\l=1$. We stress here that at present this case, along with the
particular choice $(\ref{condsuaij})$ of the coefficients
$a_{i,j}$, seems to be the only one
allowing an analytical treatment in the usual $L^p$-spaces when
dealing with a full ball.}
\end{remark}
>From $(\ref{D11})-(\ref{g22})$
we (formally) deduce that our data must
satisfy the following consistency conditions, respectively:
\begin{align}\label{DD1}
&(\textrm{C1,D})\quad\q\quad {u_0}(x)=u_1(0,x),\qq\,
&\forall&  x\in \partial\mbox{}\Omega\,,\qq\;\;\\[2mm]
\label{NN1}&(\textrm{C1,N})\q\q\q
\frac{\partial u_0}{\partial\mbox{}n}(x)=\frac{\partial u_1}{\partial\mbox{}n}
(0,x), &\forall& x\in \partial\mbox{}\Omega\,,
\qquad\quad\qquad\quad\,\\[4mm]
\label{1.18}&\hskip 2,5truecm\Phi[u_0](r)=g_1(0,r),
 &\forall&
 r\in (0,R)\,,\\[1,5mm]
\label{1.19}&\hskip 2,5truecm\Psi[u_0]=g_2(0)\,.& &\;
\end{align}
\section{Main results}
\setcounter{equation}{0}
In this section we state our {\it{local in time}}  existence and
uniqueness result related to the identification problem $\PII$.
For this purpose
we assume that the coefficients of operator
$\mathcal{A}$ satisfies
$(\ref{condsuaij})-(\ref{bcd})$, whereas, as far as the
coefficients $b_{i,j}$ and $c_i$
of operators $\cal{B},\,\cal{C}$ are concerned,
we assume:
\begin{eqnarray}\label{ipotesibijeci}
 b_{i,j}\in W^{1,\infty}(\Om)\,,
\qq c_{i}\in L^{\infty}(\Omega)\,,
\;\quad &i,j=1,2,3.&
\end{eqnarray}
In order to find out the right hypotheses on the linear operators
$\Phi$ and
$\Psi$, it will be convenient to rewrite the operator $\mathcal{A}$ in the
spherical co-ordinates $(r,\,\varphi,\,\theta)$.\\
As a consequence, using representation $(\ref{condsuaij})$ for the $a_{i,j}$'s,
through lengthy but easy computations, we  obtain the following polar representation
$\widetilde{\mathcal{A}}$ for the second-order differential
operator $\mathcal{A}$:
\begin{eqnarray}\label{tildeA}
{\widetilde{\mathcal{A}}}\!\!\!\!   & = &\!\!\!\!
D_r\big[{h}(r)D_r\big]+
\frac{2{h}(r)D_r }{r}+
\frac{{a}(r)+{b}(r)}{r^2\sin\!\theta}\Big[\,
{{(\sin\!\theta)}^{-1}D_{\varphi}^2}+D_{\theta}\big(\sin\!\theta
D_{\theta}\big)\Big]
\nonumber\\[1mm]
&   &\!\!\!  +\,\frac{1}{r^2\sin\!\theta}
\Big[\,{(\sin\!\theta)}^{-1}{D_{\varphi}\big[
\wtil{c}(r,\varphi,\theta)D_{\varphi}\big]}
+D_{\theta}\big(\wtil{c}(r,\varphi,\theta)\sin\!\theta D_{\theta}
\big)\Big]\,,
\end{eqnarray}
where we have set ${\wtil{c}}\mbox{}(r,\varphi,\theta)=
c\mbox{}(r\cos\!\varphi\sin\!\theta, r\sin\!\varphi\sin\!\theta,
r\cos\!\theta)\,$.\\
Before listing our requirements concerning operators $\Phi$ and $\Psi$ and the
data, we recall (cf. \cite{4}) some definitions about weighted Sobolev spaces.
Given an $n$-dimensional domain $\Omega$ the weighted Sobolev spaces
$W_{\sigma}^{k,p}({\Omega})$, $k\in\mathbb{N}$,
$p\in [1,+\infty]$, $\sigma\in\mathbb{R}$, are defined by
\begin{equation}\label{WSS}
W_{\sigma}^{k,p}({\Omega})=\Big\{f\in W_{loc}^{k,p}(\Omega \backslash \{0\})\,:\,
{\|f\|}_{W_{\sigma}^{k,p}({\Omega})}=
{\bigg(\sum_{0\leqslant |\alpha|
\leqslant k}\int_{{\Omega}}|x|^{\s}|D^{\alpha}f(x)|^p dx
\bigg)}^{\!1/p}\!<+\infty\Big\},
\end{equation}
\vskip -0,4truecm
\pn where
\vskip -0,5truecm
\begin{equation}
\alpha=({\alpha}_1,\ldots ,{\alpha}_n)\in \mathbb{N}^n\,,\qq
|\alpha|=\sum_{i=1}^{n}|{\alpha}_i|
\,,\qq D^{\alpha}=\frac{{\partial}^{|\alpha|}}{{\partial}^{^{{\alpha}_1}}x_1
\ldots {\partial}^{^{{\alpha}_n}}x_n}\,.
 \nonumber
\end{equation}
Of course, $W_{\sigma}^{k,p}({\Omega})$ turns out to be a Banach space
when endowed with the norm \newline
 $\|\cdot\|_{W_{\sigma}^{k,p}({\Omega})}$.
In particular, taking $\s=0$ in $(\ref{WSS})$ we obtain the usual Sobolev
spaces $W^{k,p}(\Om)$ whereas taking $k=0$  we obtain
the weighted $L^p$-spaces defined by
\begin{equation}\label{WLpS}
L_{\s}^{p}(\Om)=\Big\{f\in L_{loc}^p(\Om):\|f\|_{L_{\s}^{p}(\Om)}=\Big(
\int_{{\Omega}}|x|^{\s}|f(x)|^p dx
\Big)
^{\!1/p}\!<+\infty\Big\}.
\end{equation}

\begin{lemma}\label{suphi}
Operator $\Phi$ defined by $(\ref{Phi1})$
 maps $W^{2,p}({\Omega})$ continuously into
$W_{2}^{2,p}(0,R)$.
\end{lemma}
\begin{proof}
Taking $u\in W^{2,p}({\Omega})$ from $(\ref{Phi1})$
it follows that
\begin{equation}\label{Drj}
D_r^{(j)}\Phi[u](r)=\Phi[D_r^{(j)}u](r),\qq\forall j=0,1,2.
\end{equation}
Hence, denoting with $p'$ the conjugate exponent of $p$,
 from H\"older's inequality we obtain
\begin{align}
&\|\Phi[u]\|_{L_{2}^p(0,R)}^p=\int_0^Rr^{2}\,|
\Phi[u](r)|^p dr = \int_0^Rr^{2}\Big|\int_{\!0}^{\pi}
\!\!\sin\!\theta\mbox{}d\theta\!\int_{\!0}^{2\pi}\!\!\!u(rx')
d\varphi\Big|^p dr &\nonumber\\[2mm]
\label{0.1}&\;\,\q\qq\qq\leqslant {(4\pi)}^{p/{p'}}\!
\int_0^R\!\!r^2dr\!\!
\int_{\!0}^{\pi}\!\!\sin\!\theta d\theta\!\!
\int_{\!0}^{2\pi}\!\!\!{|u(rx')|}^p\, d\varphi
={(4\pi)}^{p/{p'}}{\|u\|}_{L^p(\Om)}^p\,.
\end{align}
Repeating similar computations and using the well-known inequalities
\begin{align}
&|D_ru(rx')|\leqslant|\nabla u(rx')|\,,\qq
|D_r^2u(rx')|\leqslant \sum_{j,k=1}^{3}|D_{x_j}D_{x_{k}}u(rx')|^2\,,&
\end{align}
from $(\ref{Drj})$ we can easily find that the following
inequalities hold:
\begin{align}
\label{1.1}&\qq\|D_r\Phi[u]\|_{L_{2}^p(0,R)}^p\leqslant
C_1{\|u\|}_{W^{1,p}(\Om)}^p\,,\qq\|D_r^2\Phi[u]\|_{L_{2}^p(0,R)}^p\leqslant
C_2{\|u\|}_{W^{2,p}(\Om)}^p\,,&
\end{align}
where $C_1$ and $C_2$ are two non-negative constants depending on  $p$
 only.\\
Therefore, from $(\ref{0.1})$ and $(\ref{1.1})$ it follows that
there exists a non-negative constant $C_3$, independent of $u$,
such that
\begin{equation}
\|\Phi[u]\|_{W_{2}^{2,p}(0,R)}\leqslant
C_3{\|u\|}_{W^{2,p}(\Om)}.
\end{equation}
\end{proof}
In this paper we will use  Sobolev spaces
 $W^{k,p}(\Omega)$ with
\begin{equation}\label{p}
p\in(3,+\infty)
\end{equation}
and we will assume that the functionals $\Phi$ and $\Psi$
satisfy the following requirements:
\begin{alignat}{5}
\label{primasuPhiePsi}
& \Phi\in\mathcal{L}\big(L^p(\Omega);\,L_{2}^p(0,R)\big),
\;\,\qq\qq
 \Psi\in L^p(\Omega)^*,&\\[1,9mm]
\label{secondasuPhi}
& \Phi[wu]=w\,\Phi[u],\qq\qq\qq\qq\;
 \forall\,(w,u)\in L_{2}^p(0,R)\times L^p(\Omega),&\\[1,9mm]
\label{terzasuPhi}
& D_r\Phi[u](r)=\Phi[D_ru](r),
\qq\qq\q \forall\,u\in W^{1,p}(\Omega)\;\;\textrm{and}\,
\,r\in(0,\,R),&
\\[1,7mm]
\label{quartasuPhi}
& \Phi\mathcal{\widetilde{A}}=\mathcal{\widetilde{A}}_1\Phi
\,\,\quad\quad\textrm{on}\;
W^{2,p}(\Omega),&\\[1,9mm]
\label{primasuPsi}
& \Psi\mathcal{\widetilde{A}}={\Psi}_1\q\quad\quad\textrm{on}\;
W^{2,p}(\Omega),\q\q\q  {\Psi}_1\in
W^{1,p}({\Omega})^*,&
\end{alignat}
where
\begin{equation}\label{tildeA1}
\mathcal{\widetilde{A}}_1= D_r\big[{h}(r)D_r]+2\frac{h(r)}{r}D_r\;.
\end{equation}
To state our result concerning the identification problem
$\PI$, we need to make also the following
assumptions on the data $f,\,u_0,\,u_1,\,g_1,\,g_2$:
\begin{alignat}{8}
\label{richiestasuf}
 &f\in C^{1+\beta}\big([0,T];L^{p}(\Omega)\big)\,,\q f(0,\cdot)\in W^{2,p}(\Om)\,,
 &\\[2,3mm]
\label{richiesteperu0eu1}
&u_0\in W^{4,p}(\Omega)\;,\q
{\cal B}u_0\in W_{\,\textrm{K}}^{2\delta,p}(\Omega)\,,\\[2,3mm]
&u_1\in C^{2+\beta}\big([0,T];L^{p}(\Omega)\big)\cap
C^{1+\beta}\big([0,T];W^{2,p}(\Omega)\big)\,,\,&\\[2,3mm]
\label{richiestaperAu0}
 &\mathcal{A}u_0+f(0,\cdot)-D_tu_1(0,\cdot)
\in W_{\textrm{K}}^{2,p}(\Omega)\,,&\\[2,3mm]
\label{richiestaperA2u0}
&F:=k_0'{\cal C}u_0+k_0{\cal B}u_0+{\mathcal{A}}^2u_0+\mathcal{A}f(0,\cdot)-D_t^2u_1(0,\cdot)+D_tf(0,\cdot)\in
W_{\,\textrm{K}}^{2\beta,p}(\Omega)\,, &
\\[2,3mm]
&g_1\in C^{2+\beta}\big([0,T];L_{2}^{p}(0,R)\big)\cap
C^{1+\beta}\big([0,T];W_{2}^{2,p}(0,R)\big),\q\; \frac{1}{r}D_tD_rg_1
C^{\beta}\big([0,T];L_{2}^{p}(0,R)\big),&\nonumber\\
\label{richiesteperg1}& &\\[2,3mm]
\label{richiesteperg2}
& g_2 \in  C^{2+\beta}\big([0,T];\mathbb{R}\big)\,,&
\end{alignat}
where $\beta\!\in\! (0,1/2)\backslash \{1/(2p)\}$, $\delta\!\in \!(\b,1/2)\backslash \{1/(2p)\}$
and function $k_0$ in $(\ref{richiestaperA2u0})$ is defined by formula $(\ref{k01})$. Moreover,
the spaces $W_{\textrm{K}}^{2,p}(\Omega)$ are defined by
\begin{equation}\label{WIK}
 W_{\textrm{K}}^{2,p}(\Omega)=\big\{
w \in  W^{2,p}(\Omega)\!: w\;
\textrm{satisfies the homogeneous condition (K)}\}\,,
\end{equation}
whereas the spaces $W_{\,\textrm{K}}^{2\g,p}(\Omega)\!\!\equiv
\!\!{\big(L^p(\Omega),
W_{\textrm{K}}^{2,p}(\Omega)\big)}_{\g,p}$,
$\g\in (0,1/2]\backslash \{1/(2p)\}$,
are interpolation
spaces between
$W_{\textrm{K}}^{2,p}(\Omega)$ and $L^p(\Omega)$ and
they are defined \cite[section 4.3.3]{5}, respectively, by:
\begin{align}
\label{WDD} & W_{\textrm{D}}^{2\g,p}(\Omega)=
\left\{\!\!\begin{array}{lll}
W^{2\g,p}(\Omega)\,, & &\,\,\,\,\;\textrm{if} \;\;
0<\g<{1}/{(2p)}\;,
\\[2mm]
\{u\in W^{2\g,p}(\Omega):u=0 \;\;\textrm{on}\;
\partial\mbox{}\Omega\}\,,& &\,\,\,\,\;\textrm{if}\;\;{1}/{(2p)}<\g\le 1/2\;,
\end{array}\right.\, &\\[3mm]
\label{WNN} & W_{\textrm{N}}^{2\g,p}(\Omega)=
W^{2\g,p}(\Omega)\,,\qq\qq\qq\qq\qq\qq\;\;\,
\textrm{if}\;\;0<\g\le 1/2\;. &
\end{align}
\begin{remark}\label{data}
\emph{Assumption $(\ref{richiesteperg1})$ ensures that
$D_t{\wtil{\cal{A}}}_1g_1\in C^{2+\b}\big([0,T], L_2^p(0,R)\big)$ (see formula
$(\ref{N10})$).}
\end{remark}
\begin{remark}
\emph{Observe that our choice $p\in(3,+\infty)$ implies the embeddings}
\begin{align}
\label{EMB1}& W^{1,p}(\Omega)\hookto C^{(p-3)/p}({\ov \Omega}),&\\
\label{EMB2}&W_2^{1,p}(0,R)\hookto C^{(p-3)/p}([0,R]).&
\end{align}
\emph{In fact, while $(\ref{EMB1})$ is a classical consequence of the Sobolev
embedding theorems (\cite{1}), Theorem 5.4, $(\ref{EMB2})$ follows immediately
from the inequalities}
\begin{eqnarray}
|u(t)-u(s)|\!\!\!&\leqslant&\!\!\! \int_s^t\!\xi^{-2/p}\xi^{2/p}|u'(\xi)|
d\xi \leqslant \bigg[\int_s^t\!\xi^{-2/(p-1)}d\xi\bigg]^{1/p'}\!\|u'\|_
{L_2^p(0,R)}\nonumber\\[2mm]
\label{hold}&\le &\!\!\!\Big(\frac{p-1}{p-3}\Big)^{\!1/p'}\!
|t-s|^{(p-3)/p}\|u\|_{W_2^{1,p}(0,R)},\qq\forall\,s,t\in[0,R]
\end{eqnarray}
\end{remark}
Assume also that $u_0$ satisfies the following conditions for some positive
constant $m$:\begin{eqnarray}
\label{J0} &J_0&\!\!\!\!\!(u_0)(r)\!:=\big|\Phi[\mathcal{C}u_0](r)\big|
\geqslant m\,,\qquad\;\forall\,r\in (0,R),\\[1,7mm]
\label{J1} &J_1&\!\!\!\!\!(u_0)\!:=\Psi[J(u_0)]\neq 0\,,
\end{eqnarray}
where we have set:
\begin{equation}\label{J}
J(u_0)(x)\!:=\!\bigg(\!\mathcal{B}u_0(x)-\frac{\Phi[\mathcal{B}u_0](|x|)}
{\Phi[\mathcal{C}u_0](|x|)}\mathcal{C}u_0(x)\!\bigg)\exp\!\bigg[
\int_{{\!|x|}}^{{R}}\!\frac{\Phi[\mathcal{B}u_0](\xi)}{\Phi[\mathcal{C}u_0]
(\xi)}d\xi\bigg]\,,\quad\forall\,x\in\Omega\,.
\end{equation}
\begin{remark}\emph{According to $(\ref{primasuPhiePsi})$ and
$(\ref{secondasuPhi})$ it follows that:
\begin{equation}
\Phi\big[J(u_0)\big](r)=\exp\!\bigg[\int_{_{r}}^{{R}}\!\frac{\Phi[
\mathcal{B}u_0](\xi)}{\Phi[\mathcal{C}u_0](\xi)}d\xi\bigg]\Phi\bigg(\!\mathcal
{B}u_0- \frac{\Phi[\mathcal{B}u_0]}{\Phi[\mathcal{C}u_0]}
\mathcal{C}u_0\!\bigg)(r)=0\,,\quad\;\forall\,r\in (0,R)\,.
\end{equation}
This means that operator $\Psi$} cannot be chosen of the form $\Psi\!=\!
\Lambda\Phi$, where\, $\Lambda$ is in $ L_{2}^p(0,R)^*$, \emph{i.e.
$\Lambda[v]\!=\!\int_0^R r^2\rho(r)v(r)dr$ for any $v\in L_{2}^p(0,R)$
and some
$\rho \in L_{2}^{p'}(0,R)$, otherwise condition $(\ref{J1})$ would be
not satisfied.
In the explicit case, when $\Phi$ and $\Psi$ have the integral representation
$(\ref{Phi1})$ and $(\ref{Psi1})$, this means that no function ${\psi}$ of
the form \,$\psi(x)=|x|^2\rho(|x|)$ is
allowed.}\end{remark}
\begin{remark}\emph{When operators $\Phi$ and $\Psi$ are defined by
$(\ref{Phi1}),\,(\ref{Psi1})$ conditions $(\ref{J0})$, $(\ref{J1})$ can be
rewritten as}:
\begin{eqnarray}
\Big| \int_0^{\pi}\!\!\!\!\sin\!\theta d\theta\!\!\int_0^{2\pi}
\!\!\!
\mathcal{C}u_0(rx') d\varphi\,\Big|\!\geqslant m_1\,,
 \qquad\forall\,r\in (0,R)\,,\qquad\qq\\
[3mm]
\qq\bigg|\int_{\!{0}}^{{R}}\!\!r^2 dr\!\!\int_0^{\pi}\!\!
\sin\!\theta d\theta\!\!
\int_0^{2\pi}\!\!\psi(rx')\bigg(\!\mathcal{B}u_0(rx')-\frac{\int_0^{\pi}
\sin\!\theta d\theta\int_0^{2\pi}\mathcal{B}u_0(rx')
d\varphi}{\int_0^{\pi}\sin\!\theta d\theta\int_0^{2\pi}
\mathcal{C}u_0(rx') d\varphi}\mathcal{C}u_0(rx')\!\bigg)\;
\nonumber\\[2,3mm]
\label{bo}\times\exp\!\Bigg[\!\int_{{r}}^{{R}}\!\!\!\frac{\;\int_0^{\pi}
\sin\!\theta d\theta\int_0^{2\pi}\mathcal{B}u_0(\xi x')
d\varphi}{\int_0^{\pi}\sin\!\theta d\theta\int_0^{2\pi}
\mathcal{C}u_0(\xi x') d\varphi}d\xi\Bigg]
d\varphi \,\bigg|\geqslant m_2\qquad\quad
\end{eqnarray}
\emph{for some positive constants $m_1$ and $m_2$.}\end{remark}
Finally, we introduce the Banach spaces ${\mathcal{U}}^{\,s,p}(T)$,
${\mathcal{U}}_{\textrm{K}}^{\,s,p}(T)$ $(\textrm{K}\in\{\textrm{D,N}\})$
which are defined for any $s\in \mathbb{N}\backslash\{0\}$ by:
\begin{equation}\label{Us}
\left\{\!\!\begin{array}{l}
{\mathcal{U}}^{\,s,p}(T)=C^s\big([0,T];L^p(\Omega)\big)\cap
C^{s-1}\big([0,T];W^{2,p}(\Omega)\big)\,,\\[2mm]
{\mathcal{U}}_{\textrm{K}}^{\,s,p}(T)=C^s\big([0,T];L^p(\Omega)\big)\cap
C^{s-1}\big([0,T];W_{\textrm{K}}^{2,p}(\Omega)\big)\,.
\end{array}\right.
\end{equation}
Moreover, we list some further consistency conditions:
\begin{align}
\label{DDV}
&\qq(\textrm{C2,D})\q\qq\q\,
{v_0}(x)=0,\hskip 2,55truecm \forall\, x\in
\partial\mbox{}\Omega\,,&\\[2mm]
\label{NNV}
&\qq(\textrm{C2,N})\q\qq\;\;\,\;\frac{\partial
v_0}{\partial\mbox{}\nu}(x)=0,
\hskip 2,37truecm\forall\, x\in
\partial\mbox{}\Omega\,,&
\end{align}
\begin{eqnarray}
\label{PHIV1}&& \Phi[v_0](r)=D_tg_1(0,r)-\Phi[D_tu_1(0,\cdot)](r),\qq\forall r\in (R_1,R_2),\\[1,5mm]
\label{PSIV1}&& \Psi[v_0]=D_tg_2(0)-\Psi[D_tu_1(0,\cdot)]\,,
\end{eqnarray}
where \begin{equation}\label{v0}
v_0(x):\!={\mathcal{A}}u_0(x)+f(0,x)-D_tu_1(0,x)\,,\qquad
\forall\,x\in\Omega\,.
\end{equation}
\begin{theorem}\label{sfera}
Let the coefficients $a_{i,j}$ $(i,j=1,2,3)$ be represented by
$(\ref{condsuaij})$ where the functions $a, b, c, d$
satisfy $(\ref{regular})$, $(\ref{bcd})$. Moreover,
let assumptions $(\ref{ipotesibijeci})$,
$(\ref{p})- (\ref{primasuPsi})$ be fulfilled and assume
that the data enjoy properties $(\ref{richiestasuf})-
(\ref{richiesteperg2})$ and satisfy
$(\ref{J0})$, $(\ref{J1})$ and the consistency conditions
$(\emph{C}1,\emph{K})$ $($cf. $(\ref{DD1})$, $(\ref{NN1}))$,
$(\emph{C}2,\emph{K})$ as well as $(\ref{1.18})$,
 $(\ref{1.19})$,  $(\ref{PHIV1})$,  $(\ref{PSIV1})$.\\
Then there exists $T^{\ast}\in (0,T]$
such that the identification problem
$\EPI$  admits a
 unique solution $(u,k)\in{\mathcal{U}}^
{\,2,p}(T^{\ast})\times C^{\beta}\big([0,T^{\ast}],
W_{2}^{1,p}(0,R)\big)$  depending continuously on the
 data with respect to the norms pointed out in
$(\ref{richiestasuf})\!
-\!(\ref{richiesteperg2})$.\\
In the case of the specific operators $\Phi$, $\Psi$ defined by $(\ref{Phi1}),
\,(\ref{Psi1})$ the previous results are still true if
   $\psi\in
C^1(\overline{\Omega})$, with
 $\psi_{|_{\partial{\Omega}}}\!=\!0$ when
$\emph{K}\!=\!\emph{D}$.
\end{theorem}
\begin{corollary}\label{PHIPSIBAll}
 When $\Phi$ and $\Psi$ are  defined by $(\ref{Phi1})$ and
$(\ref{Psi1})$,
respectively, and the coefficients $a_{i,j}\;(i,j=1,2,3)$ are
 represented by $(\ref{condsuaij})$, conditions
$(\ref{primasuPhiePsi})-(\ref{primasuPsi})$ are
 satisfied under assumptions $(\ref{regular})$, $(\ref{p})$
 and the hypothesis
$\psi\in C^1(\overline{\Omega})$, with
 ${\psi}_{|_{\partial\mbox{}\Omega}}\!=\!0$ when
$\emph{K}\!=\!\emph{D}$.
\end{corollary}
\begin{proof}
>From definitions $(\ref{Psi1})$ and H\"older's inequality it
immediately follows
\begin{equation}\label{psinorm}
\big|\Psi[v]\big|\leqslant
{\|\psi\|}_{C(\overline{\Om})} {\|v\|}_{L^1(\Om)}\leqslant
{\bigg[\frac{4}{3}\pi R^3\bigg]}^{\!1/{p'}}\!
{\|\psi\|}_{C(\overline{\Om})}{\|v\|}_{L^p(\Om)}\,.\qq\qq
\end{equation}
Hence, from $(\ref{0.1})$ and $(\ref{psinorm})$ we have that
$(\ref{primasuPhiePsi})$ is satisfied. Definition
$(\ref{Phi1})$ easily implies $(\ref{secondasuPhi})$ and
$(\ref{terzasuPhi})$, as we have already noted in $(\ref{Drj})$.
So, it remains only to prove
that decompositions $(\ref{quartasuPhi})$ and
$(\ref{primasuPsi})$ hold.\\
When the coefficients $a_{i,j}$ are represented by
$(\ref{condsuaij})$ the second-order differential operator
$\mathcal{A}$ can be represented, in spherical co-ordinates,
by operator $\widetilde{\mathcal{A}}$
defined by $(\ref{tildeA})$.
Our next task consists in computing $\Phi\big[\wtil{\cal{A}}w\big]$
for any $w\in W_{\textrm{K}}^{2,p}(\Om)$, $p\in(3,+\infty)$.
Observe first that from $(\ref{0.1})$ and $(\ref{tildeA1})$ it follows
\begin{equation}\label{primo}
\Phi[\widetilde{\mathcal{A}}_1w](r)
 =   \int_0^{\pi}\!\!\!\!\sin\!\theta d\theta\!\!\int_0^{2\pi}\!\!\!
\lambda({Rx}')\Big\{\! D_r[ h(r)D_r w(rx')]
+2\frac{h(r)}{r}D_rw(rx')\!\Big\} d\varphi=
\widetilde{\mathcal{A}}_1\Phi[w](r)
\end{equation}
\vskip -0,3truecm
\pn Since $p\in (3,+\infty)$, using the Sobolev embedding theorem of $W^{1,p}(\Om)$ into $C(\ov\Om)$
 and the well-known formulae
\begin{equation}
\label{Dr}
\left\{\!\! \begin{array}{lll}
D_{r}\!\!\! & = &\!\!\!{\cos\!\varphi\sin\!\theta} D_{x_1}+\sin\!\varphi
\sin\!\theta D_{x_2}+\cos\!\theta D_{x_3}\,,\\[1,7mm]
D_{\varphi}\!\!\!  & = &\!\!\!{-r\sin\!\varphi\sin\!\theta}  D_{x_1}+
r\cos\!\varphi\sin\!\theta D_{x_2}\,, \\[1,7mm]
D_{\theta} \!\!\! & = &\!\!\!{r\cos\!\varphi\cos\!\theta}
D_{x_1}+r\sin\!\varphi\sin\!\theta D_{x_2}-r\sin\!\theta D_{x_3}\,,
\end{array}\right.
\end{equation}
it can be easily shown that  $(D_{\varphi}w)/(r\sin\!\theta)$
and $(D_{\theta}w)/r$ are bounded, while the functions $(D_{\varphi}^2w)/\sin\!\theta$
and $D_{\theta}(\sin\!\theta D_{\theta}w)$ belong to
$L^1(\partial\mbox{} B(0,r))$ for every $r\in (0,R)$.
Therefore, integrating by parts, we obtain
\begin{align}
&\Phi\Big[\frac{{a}(r)+{b}(r)}{r^2\sin\!\theta}\Big(
{(\sin\!\theta)}^{-1}{D_{\varphi}^2w}+D_{\theta}(\sin\!\theta
D_{\theta}w)\!\Big)\!\Big]\!(r)\nonumber\\[3mm]
\label{terzo}&\q =\frac{{a}(r)+{b}(r)}{r^2}\bigg\{
\!\int_0^{\pi}\!\bigg[\frac{D_{\varphi}w(rx')}{\sin\!\theta}
\bigg|_{\varphi=0}^{\varphi=2\pi}\,\bigg] d\theta+\int_{0}^{2\pi}\!\Big[
{D_{\theta}w(rx')\sin\!\theta]}\Big|_{\theta=0}^{\theta=\pi}\,\Big]
d\varphi\bigg\}=0
\,,\\[3mm]
&\Phi\Big[\,\frac{1}{r^2\sin\!\theta}
\Big({(\sin\!\theta)}^{-1}{D_{\varphi}\big[
\wtil{c}(r,\varphi,\theta)D_{\varphi}w\big]}
+D_{\theta}\big[\wtil{c}(r,\varphi,\theta)\sin\!\theta D_{\theta}w
\big]\Big)\!\Big]\!(r)\nonumber\\[3mm]
\label{quarto}&\q =\frac{1}{r^2}\bigg\{
\!\int_0^{\pi}\!\bigg[\frac{\wtil{c}(r,\varphi,\theta)D_{\varphi}w(rx')}
{\sin\!\theta}\bigg|_{\varphi=0}^{\varphi=2\pi}\,\bigg] d\theta+
\int_{0}^{2\pi}\!\Big[
{\wtil{c}(r,\varphi,\theta)D_{\theta}w(rx')\sin\!\theta]}
\Big|_{{\theta=0}}^{{\theta=\pi}}\,\Big] d\varphi\bigg\}=0.
\end{align}
Hence, from  $(\ref{primo})$, $(\ref{terzo})$,
$(\ref{quarto})$ we find
that $(\ref{quartasuPhi})$  holds for every
$w\in W_{\textrm{K}}^{2,p}(\Omega)$
with $p\in (3,+\infty)$.\\
Let now  $\Psi$ be the functional defined in
$(\ref{Psi1})$. Analogously to what we have done for $\Phi$, we apply $\Psi$
to both sides in  $(\ref{tildeA})$. Performing computations similar to those
made above and using the assumption
$\psi_{|_{\partial\mbox{}\Omega}}\!=\!0$ when $\textrm{K}=\textrm{D}$
which ensure that the surface integral vanishes,
we obtain the equation
$$\Psi[\widetilde{\cal{A}}w]= {\Psi}_1[w]\,,\qq w\in W_{\textrm{K}}^{2,p}(\Omega)\,,$$
\vskip -0,3truecm
\pn where
\vskip -0,3truecm
\begin{align}
& {\Psi}_1[ w]=
 -\!\int_{0}^{R}\!\!r^2\,{h}(r)dr\!\!
\int_0^{\pi}\!\!\!\!\sin\!\theta d\theta\!\!
\int_0^{2\pi}\!\!\!D_rw(rx')D_r\psi(rx')d\varphi&\nonumber
\\[1,7mm]
&\qq\q-\!\int_{0}^{R}\!\!\!r^2dr\!\!\int_0^{\pi}\!\!\!\!\sin\!\theta
d\theta\!\!
\int_0^{2\pi}\!\big[{a}(r)+{b}(r)+\widetilde{c}(r,\varphi,\theta)
\big]\frac{D_{\varphi}w(rx')}
{r\sin\!\theta} \frac{D_{\varphi}\psi(rx')}{r\sin\!\theta}
\,d\varphi&\nonumber
\\[1,7mm]
\label{psi11ball}&\qq\q -\int_{0}^{R}\!\!\!r^2dr\!\!
\int_0^{\pi}\!\!\!\!\sin\!\theta d\theta\!\!\int_0^{2\pi}\!\big[{a}(r)
+{b}(r)+\widetilde{c}(r,\varphi,\theta)
\big]\frac{{D_{\theta}w(rx')}}{r}
\frac{D_{\theta}\psi(rx')}{r}\,d\varphi\,.&
\end{align}
Now it is an easy task to show that $\Psi_1$ defined in $(\ref{psi11ball})$
belongs to ${W^{1,p}(\Om)}^{\ast}$. Indeed, using  formulae
$(\ref{Dr})$ and H\"older's inequality, we can easily find
\begin{equation}\label{C1}
|\Psi_1[w]|\leqslant
{C_1}\|\nabla u\|_{L^{p}(\Om)}\leqslant C_1\|w\|_{W^{1,p}(\Om)}\,,
\end{equation}
where $C_1>0$ depends
 on ${\|\psi\|}_{C^{1}(\ov\Om)}$ and
$\max\!\big[\|h\|_{L^{\infty}(0,R)},{\|a+b+c\|}_{L^{\infty}(\Om)}\big]$,
only.\\
Hence also decomposition $(\ref{primasuPsi})$ holds and this
completes the proof.\ \end{proof}
\section{An equivalence result in the concrete case}
\setcounter{equation}{0}
Taking advantage of the results proved in \cite{2}, we limit ourselves
to sketching the procedure for solving the necessary equivalence result.\\
We introduce the new triplet of unknown functions $(v, l, q)$ defined by
\begin{eqnarray}\label{v,h,q}
 v(t,x)= D_tu(t,x)-D_t{u}_1(t,x)\,,\q\; l(t)=k(t,R_2)\,,\q\;
 q(t,r)=D_rk(t,r)\,,\q
\end{eqnarray}
so that  $u$ and $k$ are given, respectively, by the following formulae
\begin{eqnarray}
 u(t,x)\!\!\!&=&\!\!\!u_1(t,x)-u_1(0,x)+u_0(x) +\int_0^t\!v(s,x)ds,\q\;
\forall\,(t,x)\in[0,T]\times\Om,\q\\[1mm]
k(t,r)\!\!\!&=&\!\!\!l(t)-\int_{r}^{R}\!\!\!\!q(t,\xi)d\xi:=l(t)-Eq(t,r),\quad
\,\forall\,
(t,r)\in[0,T]\times(0,R).
\end{eqnarray}
Then problem $(\ref{problem})$, $(\ref{u0})-(\ref{g22})$ can be shown
to be equivalent to the following identification problem:
\begin{eqnarray}\label{problem1}
D_tv(t,x)\!\!\!&=&\!\!\!\mathcal{A}v(t,x)+\int_0^t\! k(t-s,|x|)
\big[\mathcal{B}v(s,x)+\mathcal{B}D_{t}u_1(s,x)\big]ds+k(t,|x|)
\mathcal{B}u_0(x)\nonumber\\[1,2mm]
\!\!\!& &\!\!\!+\!\int_0^t\! D_{|x|}k(t-s,|x|)\big[\mathcal{C}v(s,x)+
\mathcal{C}D_t{u}_1(s,x)\big]ds+D_{|x|}k(t,|x|)\mathcal{C}u_0(x)
\nonumber\\[1,7mm]
&&\!\!\!+ \mathcal{A}D_t{u}_1(t,x)-D_t^{2}{u}_1(t,x)
+D_tf(t,x),\qq\forall\,(t,x)\in[0,T]\times\Omega,\q\;\;\q
\end{eqnarray}
\vskip -0,85truecm
\begin{align}
\label{v01}&\q v(0,x)={\mathcal{A}}u_0(x)+f(0,\cdot)-D_tu_1(0,x)\!:=v_0(x),\qq
\forall\,x\in\Omega,&\\[1,5mm]
 &\q v\; \textrm{satisfies the homogeneous boundary condition (K)},
\q\textrm{K}\in\{\textrm{D,N}\}\,,&\\[3mm]
\label{hhh}
&\q l(t)= l_0(t)+N_3(v,l,q)(t),\,\qq \forall\;t\in[0,T],
 \\[3mm]
\label{q3}&\q q(t,r)=q_0(t,r)+J_2(u_0)(r)N_3(v,l,q)(t)
+N_2(v,l,q)(t,r),\q\forall\,(t,r)\in [0,T]\times (0,R),\nonumber\\
\end{align}
\vskip -0,3truecm
\pn where we have set
\begin{align}
\label{h0}
& l_0(t)\!:={[J_1(u_0)]}^{-1}N_0(u_0,u_1,g_1,g_2,f)(t)\,,
\qq \forall\;t\in[0,T],\\[2,3mm]
& \label{q0}
q_0(t,r)\!:=J_2(u_0)(r)h_0(t)+N_3^0(u_0,u_1,g_1,f)(t,r),\;\quad\forall\,(t,r)
\in [0,T]\times (0,R).
\end{align}
We recall that operators $J_0$, $J_1$ and $J_2$ are defined, respectively,
by $(\ref{J0})$, $(\ref{J1})$ and
\begin{equation}\label{J2}
 J_2(u_0)(r)=-\frac{\Phi[\mathcal{B}u_0](r)}
{\Phi[\mathcal{C}u_0](r)}\exp\!
\bigg[\!\int_{\!{r}}^{\!{R_2}}\frac{\Phi[\mathcal{B}u_0](\xi)}{\Phi[
\mathcal{C}u_0](\xi)}d\xi\bigg],\qquad\forall\,r\in (0,R).
\end{equation}
To define operators $N_2$ and $N_3$ appearing in $(\ref{hhh})$,
$(\ref{q3})$ we need to introduce the operators $N_1$ and $L$:
\begin{align}
&\; {N}_1(v,l,q)(t,|x|):=-\!\int_0^t\!
\big[ l(t-s)-Eq(t-s,|x|)
\big]\big[\mathcal{B}v(s,x)+\mathcal{B}D_{t}u_1(s,x)\big]ds
&\nonumber\\
\label{N1}&\qq\qq -\int_0^t\!\!q(t-s,| x |)\big
[\mathcal{C}v(s,x)+\mathcal{C}
D_t{u}_1(s,x)\big]ds\,,\;\quad\forall\;(t,x)\in[0,T]\!\times\!\Omega\,,&\\[3mm]
\label{L}
&\, Lg(t,r)\!:=\int_{r}^{R_2}\!\!\!\exp\!\bigg[\int_{\!r}^{\eta}\frac{\Phi[
\mathcal{B}u_0](\xi)}{\Phi[\mathcal{C}u_0](\xi)}d\xi\bigg]\frac{g(t,\eta)}
{\Phi[\mathcal{C}u_0](\eta)}d\eta\,,\q\;\forall g\in L^1((0,T)\times (0,R)).&
\end{align}
Now, denoting by $I$ the identity operator, define $N_2$ and $N_3$
via the formulae
\begin{eqnarray}
N_2(v,l,q)(t,r)\!:\!\!\!\!\!\!&=&\!\!\!\frac{1}
{\Phi[\mathcal{C}u_0](r)}\big[I+\Phi[
\mathcal{B}u_0](r)L\big]\,\Phi[{N_{1}}(v,l,q)(t,\cdot)](r)
\nonumber\\[2mm]
\label{N2} &: =&\!\!\!\!J_3(u_0)(r)
\,\Phi[{N_{1}}(v,l,q)(t,\cdot)](r)
,\qquad\qquad\\[3mm]
 N_3(v,l,q)(t)\!:\!\!\!\!\!\!&=&\!\!\!{[J_1(u_0)]}^{-1}
\Big\{\Psi[{N_{1}}(v,l,q)(t,\cdot)]\!-\!\Psi[N_2(v,l,q)(t,\cdot)
\mathcal{C}u_0]\nonumber\\[2mm]
\label{N3} & +&\!\!\!\!\Psi\big[E\big(N_2(v,l,q)(t,\cdot)\big)
\mathcal{B}u_0\big]\!-\!{\Psi}_1[v(t,\cdot)]\Big\}\,,
\end{eqnarray}
where $\Psi_1$ is defined by  $(\ref{psi11ball})$.\\
Finally, to define operators $N_0$ and $N_3^0$  appearing in
$(\ref{h0})$, $(\ref{q0})$ we need to introduce first the operators $N_1^0$
and $N_2^0$, where operators ${\wtil{\cal{A}}}$ and  ${\wtil{\cal{A}}}_1$
are defined, respectively, by $(\ref{tildeA})$ and $(\ref{tildeA1})$:
\begin{eqnarray}
\hskip -0,7truecm N_1^{0}(u_1,g_1,f)(t,r)\!\!\!&=&\!\!\!D_t^2g_1(t,r)
-D_t{\widetilde{\mathcal{A}}}_1g_1(t,r)
\nonumber\\[2mm]
\hskip -0,7truecm &&\!\!\!-\Phi[D_tf(t,\cdot)](r)\,,\qq\forall\,(t,r)\in[0,T]\!
\times\!(0,R),\label{N10}\\[3,5mm]
\label{N20}\hskip -0,7truecm N_2^{0}(u_1,g_2,f)(t)\!\!\!&
=&\!\!\!\!D_t^2g_2(t)-{\Psi}_1[D_tu_1(t,\cdot)]
-{\Psi}[D_t f(t,\cdot)]\,,\qq\forall\,t\in[0,T]\,.
\end{eqnarray}
Then we define
\begin{eqnarray}
N_3^0(u_0,u_1,g_1,f)(t,r)\!\!\!\!&:=&\!\!\!
\frac{1}{\Phi[\mathcal{C}u_0](r)}
\big[I+\Phi[\mathcal{B}u_0](r)L\big]N_1^{0}(u_1,g_1,f)(t,r)
\nonumber\\[2mm]
\label{N30}&:=&\!\!\!\!J_3(u_0)(r)N_1^{0}(u_1,g_1,f)(t,r),\\[3mm]
N_0(u_0,u_1,g_1,g_2,f)(t)\!:\!\!\!&=&\!\!\!N_2^{0}(u_1,g_2,f)(t)-
\Psi[N_3^0(u_0,u_1,g_1,f)(t,\cdot)\mathcal{C}u_0]\nonumber\\[1,8mm]
\label{N0}& &\!\!\!-\Psi\big[E\big(N_3^0(u_0,u_1,g_1,f)(t,\cdot)\big)
\mathcal{B}u_0\big]\,.
\end{eqnarray}
Finally, we introduce function $k_0$ appearing in $(\ref{richiestaperA2u0})$:
\begin{eqnarray}
\label{k01} k_0(r)\!\!\!\!&=&\!\!\!\![J_1(u_0){]}^{-1}\Big\{
\Psi[\wtil{l}_2]
+N_2^{0}(u_1,g_2,f)(0)-{\Psi}_1[v_0]\!\Big\}
\exp\!\bigg[\int_{\!r}^{R_2}\frac{\Phi[
\mathcal{B}u_0](\xi)}{\Phi[\mathcal{C}u_0](\xi)}d\xi\bigg]
\nonumber\\[1,2mm]
&&\!\!\!\!+\int_{\!R_2}^{{r}}
\!\!\!\exp\!\bigg[\int_{\!r}^{\eta}\!\frac{\Phi[\mathcal{B}u_0](\xi)}
{\Phi[\mathcal{C}u_0](\xi)}d\xi\bigg]\frac{N_1^0(u_1,g_1,f)(\eta)}
{\Phi[\mathcal{C}u_0](\eta)}d\eta\,,\q \;\forall\;r\in (R_1,R_2)\,.
\end{eqnarray}
where for any $x\in\Omega$ we set
\begin{eqnarray}
\wtil{l}_2(x)\!\!\!\!&:= &\!\!\!\!\mathcal{C}u_0(x)\bigg\{
\frac{N_1^0(u_1,g_1,f)(|x|)}
{\Phi[\mathcal{C}u_0](|x|)}-\frac{\Phi[\mathcal{B}u_0](|x|)}
{\Phi[\mathcal{C}u_0](|x|)}\int_{\!R_2}^{|x|}\!\!\!\exp\!
\bigg[\int_{\!|x|}^{\eta}\frac{\Phi[\mathcal{B}u_0](\xi)}
{\Phi[\mathcal{C}u_0]
(\xi)}d\xi\bigg]\nonumber\\\nonumber\\[2mm]
\label{l2}&
&\!\!\!\!\times\frac{N_1^0(u_1,g_1,f)(\eta)}{\Phi[\mathcal{C}u_0](\eta)}
d\eta\bigg\}+\mathcal{B}u_0(x)\int_{\!R_2}^{|x|}\!\!\!\exp\!\bigg[
\int_{\!|x|}^{\eta}\frac{\Phi[\mathcal{B}u_0](\xi)}
{\Phi[\mathcal{C}u_0](\xi)}d\xi\bigg]\frac{N_1^0(u_1,g_1,f)(\eta)}
{\Phi[\mathcal{C}u_0](\eta)}d\eta\,.\nonumber
\end{eqnarray}
We can summarize the result sketched in this section in the following
equivalence theorem.
\begin{theorem}\label{3.1}
The pair  $(u,k)\in{\mathcal{U}}^{\,2,p}(T)\times
C^{\beta}\big([0,T];W_{2}^{1,p}
(0,R)\big)$ is a solution to the identification problem $\emph{\PII},\;
\emph{K}\in\{\emph{D,N}\}$, if and only if the triplet
$(v,l,q)$ defined by
$(\ref{v,h,q})$  belongs to $ {\mathcal{U}}_{\emph{K}}^{\,1,p}(T)
\times C^{\beta}
\big([0,T];\mathbb{R}\big)\times C^{\beta}\big([0,T];
L_{2}^{p}(0,R)\big)$
and solves problem $(\ref{problem1})\!-\!(\ref{q3})$.
\end{theorem}
\section{An abstract formulation of problem
(\ref{problem1})-(\ref{q3}).}
\setcounter{equation}{0}
Starting from the result of the previous section, we can reformulate our
identification problem in a Banach space framework.\\Let
$A:\mathcal{D}(A)\subset X \to X$ be a linear closed operator satisfying the
following assumptions:
\begin{itemize}
\item[(H1)]\emph{there exists $\zeta\in (\pi /2,\pi)$ such that the resolvent
set of $A$ contains $0$ and the open sector ${\Sigma}_{\zeta}=\{\mu\in
\mathbb{C}:|\arg\mu|<\zeta\}$;}
\item[(H2)]\emph{there exists $M>0$ such that ${\|{(\mu I-A)}^{-1}\|}
_{\mathcal{L}(X)}\leqslant M|\mu{|}^{-1}$ for every $\mu\in {\Sigma}_{\zeta}$;}
\item[(H3)]\emph{$X_1$ and $X_2$ are Banach spaces such that
$\mathcal{D}(A)=X_2\hookrightarrow X_1\hookrightarrow X $. Moreover,
$\mu\to {(\mu I-A)}^{-1}$ belongs to ${\cal L}(X;X_1)$ and satisfies the
estimate ${\|{(\mu I-A)}^{-1}\|}
_{\mathcal{L}(X;X_1)}\leqslant M|\mu{|}^{-1/2}$ for every $\mu\in
{\Sigma}_{\zeta}$.}
\end{itemize}
Here $\mathcal{L}(Z_1;Z_2)$ denotes, for any pair of Banach spaces
$Z_1$ and $Z_2$, the Banach space of all bounded linear operators
 from $Z_1$ into $Z_2$ equipped with the uniform-norm.
 In particular we set ${\cal L}(X)=\mathcal{L}(X;X)$.\\
By virtue of assumptions (H1), (H2) we can define the analytic semigroup
$\{{\rm e}^{tA}\}_{t\geqslant 0}$ of bounded linear operators in
$\mathcal{L}(X)$
 generated by $A$. As is well-known, there exist positive constants
$\widetilde{c_{k}}(\zeta)\; (k \in\mathbb{N})$ such that
$$
\|A^k{\rm e}^{tA}\|_{\mathcal{L}(X)}\leqslant \widetilde{c_{k}}(\zeta)Mt^{-k},
\qquad\forall t \in {\mathbb{R}}_{+},\, \forall k\in\mathbb{N}.
$$
After endowing $\mathcal{D}(A)$ with the graph-norm,
we can define the following family of interpolation spaces
${\mathcal{D}}_{A}(\beta,p)$, $\beta\in (0,1)$, $p\in [1,+\infty]$, which are
intermediate between $\mathcal{D}(A)$ and $X$:
\begin{eqnarray}\label{interpol1}
{\mathcal{D}}_{A}(\beta,p)=
\Big\{x\in X: |x|_{{\mathcal {D}}_{A}(\beta,p)} < +\infty\Big\}, \qq \mbox{if } p\in [1,+\infty],
\end{eqnarray}
where
\begin{equation}
{|x|}_{{\mathcal{D}}_{A}(\beta,p)} = \left\{
\begin{array}{l}
\ds \Big(\int_0^{+\infty}\!t^{(1-\beta)p-1}\|A{\rm e}^{tA}x\|_X^p\,dt
\Big)^{\! 1/p},\q \mbox{if } p\in [1,+\infty),
\\[5mm]
\sup_{0<t\le 1}\big(t^{1-\beta}\|A{\rm{e}}^{tA}x\|_X\big),\q \hskip 1.2truecm \mbox{if } p=\infty.
\end{array} \right.
\end{equation}
\pn
They are well defined by virtue of assumption (H1). Moreover, we set
\begin{equation}\label{interpol2}
{\mathcal{D}}_{A}(1+\beta,p)\!=\! \{x\in\mathcal{D}(A):Ax\in
{\mathcal{D}}_{A}(\beta,p)\}\,.
\end{equation}
Consequently, ${\mathcal{D}}_{A}(n+\beta,p)$, $n\in\mathbb{N},
\beta\in (0,1)$, $p\in [1,+\infty]$, turns out to be a Banach space when equipped
with the norm
\begin{equation}
{\|x\|}_{{\mathcal{D}}_{A}(n+\beta,p)}\!=\!
\sum_{j=0}^{n}{\|A^{j}x\|}_{X}+{|A^{n}x|}_{{\mathcal{D}}_{A}(\beta,p)}\,.
\end{equation}
In order to reformulate in an abstract form our identification
problem$(\ref{problem1})\!-\!(\ref{q3})$
 we need the
following assumptions involving spaces, operators and data:
\begin{alignat}{9}
&(\textrm{H}4)\;\emph{$Y$ and $Y_1$ are Banach spaces such that
$Y_1\hookrightarrow Y$;}&
\nonumber\\[2mm]
&(\textrm{H}5)\;\emph{$B:\mathcal{D}(B)\subset X\rightarrow X$ is a linear
closed operator such that $X_2\subset \mathcal{D}(B)$;}&
\nonumber\\[2mm]
&(\textrm{H}6)\;\emph{$C:\mathcal{D}(C):=X_1\subset X\rightarrow
X$ is a linear closed operator;}&
\nonumber\\[2mm]
&(\textrm{H}7)\;\emph{$E\in\mathcal{L}(Y;Y_1)$, $\Phi\in\mathcal{L}(X;Y)$,
 $\Psi\in {X}^{\ast}$, ${\Psi}_{1}\in
{X_1}^{\ast}$;}&
\nonumber\\[2mm]
&(\textrm{H}8)\;\emph{$\mathcal{M}$ is a continuous bilinear operator
from $Y\times {\wtil X}_1$ to $X$ and from $Y_1\times X$ to $X$,}&
\nonumber\\
& \qq\;\emph{where $X_1\hookto {\wtil X}_1$;}\nonumber\\[2mm]
&(\textrm{H}9)\;\emph{$J_1:X_2\rightarrow\mathbb{R}$,
$J_2:X_2\rightarrow Y$, $J_3:X_2\rightarrow\mathcal{L}(Y)$\,
are three prescribed (non-linear)}&
\nonumber\\
&\qq\;\emph{operators}\,;&\nonumber\\[1mm]
&(\textrm{H}10)\;\emph{$u_0, v_0\in X_2$,\, $Cu_0\in X_1$,\;
$J_1(u_0)\neq 0$, $Bu_0\in \mathcal{D}_A(\delta,+\infty)$, $\delta\in (\b,1/2)$ ;}&
\nonumber\\[2mm]
&(\textrm{H}11)\;\emph{$q_0\in C^{\beta}([0,T];Y)$,
$l_0\in  C^{\beta}([0,T];\mathbb{R})$\,;}
&\nonumber\\[2mm]
&(\textrm{H}12)\;\emph{$z_0\in C^{\beta}([0,T];X)$,\;
$z_1\in C^{\beta}([0,T];{\wtil X}_1)$,\; $z_2\in C^{\beta}([0,T];X)$\,;}&\nonumber\\[2mm]
&(\textrm{H}13)\;\emph{$Av_0+\mathcal{M}({\wtil q}_0,{C}u_0)+
{\wtil l}_0{B}u_0-\mathcal{M}(E{\wtil q}_0,Bu_0)+z_2(0,\cdot)
\in\mathcal{D}_A(\beta,+\infty)$\,.}&\nonumber
\end{alignat}
The elements ${\wtil q}_0$ and ${\wtil l}_0$ appearing in
$(\textrm{H}13)$ are defined by:
\begin{equation}\label{rem4.2.1}
\left\{\!\begin{array}{l}
\wtil{l}_0=l_0(0)-\big[J_1(u_0)\big]^{-1}\Psi_1[v_0]\,,\\[3mm]
\wtil{q}_0=q_0(0)+J_2(u_0)\big[J_1(u_0)\big]^{-1}\Psi_1[v_0]\,,
\end{array}\right.
\end{equation}
where $l_0$ and $q_0$ are the elements appearing in $(\textrm{H}11)$.
\begin{remark}\label{rem4.3}
\emph{In the explicit case we get the equations
\begin{equation}
\wtil{l}_0=k_0(R_2)\,,\,\q\wtil{q}_0(r)=k_0'(r)\,.
\end{equation}
where $k_0$ is defined in $(\ref{k01})$.}
\end{remark}
We can now reformulate our direct problem:
 \emph{determine a function
$v\in C^1([0,T];X)\cap C([0,T];X_2)$ such that}
\begin{eqnarray}
\label{problem2}
v'(t)\!\!\!&=&\!\!\![\l_0I+A]v(t)+\!\int_0^t\!\! l(t-s)[{B}v(s)+z_0(s)]ds-
\!\int_0^t\!\! \mathcal{M}\big(Eq(t-s),{B}v(s)+z_0(s)\big)ds\nonumber\\[2mm]
&& + \int_0^t \mathcal{M}\big(q(t-s),{C}v(s)+z_1(s)\big)ds
+\mathcal{M}\big(q(t),{C}u_0\big)+l(t)Bu_0\nonumber\\[2mm]
& &-\mathcal{M}\big(Eq(t),Bu_0\big)+z_2(t),
\hskip 3truecm \forall\;t\in[0,T],\\[2mm]
\label{v02}
v(0)\!\!\!&=&\!\!\!v_0.
\end{eqnarray}
\begin{remark}\label{z0z1z2}
\emph{In the explicit case $(\ref{problem1})-(\ref{q3})$ we have
$A={\cal{A}}-\l_0I$, with a large enough positive $\l_0$, and the
functions
 $z_0, z_1, z_2$ defined by}
\begin{eqnarray}
\label{z1z2z3}&z_0=D_t\mathcal{B}u_1\;,\qquad
z_1=D_t\mathcal{C}u_1\;,\qquad
 z_2=D_t\mathcal{A}u_1-D_t^2u_1+D_tf,&
\end{eqnarray}
\emph{whereas $v_0, h_0, q_0$ are defined, respectively,
 via the formulae $(\ref{v0})$, $(\ref{h0})$, $(\ref{q0})$.}
\end{remark}
Introducing the operators
\begin{eqnarray}
\widetilde{R}_2(\!\!\!\!\!&v&\!\!\!\!\!,h,q)\!:=-{[J_1(u_0)]}^{-1}
\Big\{\!\Psi\big[\mathcal{M}\big(J_3(u_0)\Phi[{N_{1}}(v,l,q)],Cu_0\big)\big]
\nonumber\\[1,8mm]
\label{tildeR2}&&\qquad\qquad-\Psi\big[\mathcal{M}\big(E\big(J_3(u_0)
\Phi[{N_{1}}(v,l,q)],{B}u_0\big)\big]-\Psi[{N_{1}}(v,l,q)]\!\Big\}\,,
\\[1,8mm]
\label{tildeR3}\widetilde{R}_3(\!\!\!\!\!&v&\!\!\!\!\!,h,q)\!:=J_2(u_0)
\widetilde{R}_2(v,l,q)+J_3(u_0)\Phi[{N_{1}}(v,l,q)]\,,\\[1,8mm]
\widetilde{S_2}(\!\!\!\!\!&v&\!\!\!\!\!)\!:=\!{[J_1(u_0)]}^{-1}
\Big\{\!\Psi\big[\mathcal{M}\big(J_3(u_0){\Phi}_1[v],Cu_0\big)\big]\!+\!
\Psi\big[\mathcal{M}\big(E\big(J_3(u_0){\Phi}_1[v],Cu_0\big)\big]\!
-\!{\Psi}_1[v]\!\Big\}\,,\nonumber\\\label{tildeS2} \\
\label{tildeS3}\widetilde{S_3}(\!\!\!\!\!&v&\!\!\!\!\!)\!:=J_2(u_0)
\widetilde{S_2}(v)\,,
\end{eqnarray}
the fixed-point system $(\ref{hhh})$, $(\ref{q3})$ for $l$ and $q$ becomes
\begin{eqnarray}
\label{ha1}l=l_0+\widetilde{R}_2(v,l,q)+\widetilde{S_2}(v)\,,\\[1,6mm]
\label{qa1}q=q_0+\widetilde{R}_3(v,l,q)+\widetilde{S_3}(v)\,.
\end{eqnarray}
The present situation is analogous to the one in \cite{3} (cf. Section 4).
Consequently, also in this case we can apply the abstract results
proved in \cite{2} (cf. Sections 5 and 6) to get the following local
in time existence and uniqueness theorem.
\begin{theorem}\label{4.2}
Under assumptions $(\emph{H}1)-(\emph{H}13)$
there exists $T^{\ast}\in (0,T)$
such that for
any $\tau\in (0,T^{\ast}]$ problem $(\ref{problem2}), (\ref{v02}),
(\ref{ha1}), (\ref{qa1})$ admits a unique solution $(v,l,q)\in
[C^{1+\beta}([0,\tau];X)\cap C^{\beta}([0,\tau];X_2)]\times C^{\beta}([0,\tau];\mathbb{R})\times
C^{\beta}([0,\tau];Y)$.
\end{theorem}
\section{Solving the identification problem (\ref{problem1})--(\ref{q3})\newline and proving Theorem \ref{sfera}}
\setcounter{equation}{0}
The main difficulties we meet when we try to solve our identification problem
$\PI$, in the open ball $\Omega$ can be overcome by introducing the
representation $(\ref{condsuaij})$ and  the additional assumptions
$(\ref{regular})\!-\!(\ref{bcd})$ for the coefficients
$a_{i,j}\;(i,j=1,2,3)$ of $\mathcal{A}$.

The basic result of this section is the following Theorem.
\begin{theorem}\label{sfera1}
Let the coefficients $a_{i,j}$ $(i,j=1,2,3)$ be represented by
$(\ref{condsuaij})$ where the functions $a, b, c, d$
satisfy $(\ref{regular})\!-\!(\ref{bcd})$. Moreover,
let assumptions $(\ref{ipotesibijeci})$,
$(\ref{p})\!-\!(\ref{richiesteperg2})$,
$(\ref{J0})$, $(\ref{J1})$ be fulfilled along with the
consistency conditions $(\ref{DDV})-(\ref{PSIV1})$.
\\
Then there exists $T^{\ast}\in (0,T]$
such that the identification problem
$(\ref{problem1})-(\ref{q3})$  admits a
 unique solution $(v,l,q)\in{\mathcal{U}}_{\emph{K}}^
{\,1,p}(T^{\ast})\times C^{\beta}\big([0,T^{\ast}];\mathbb{R}\big)\times C^{\beta}\big([0,T^{\ast}];
L_{2}^{p}(0,R)\big)$  depending continuously on the
 data with respect to the norms pointed out in
$(\ref{richiestasuf})\!
-\!(\ref{richiesteperg2})$.\\
In the case of the specific operators $\Phi$, $\Psi$ defined by $(\ref{Phi1}),
\,(\ref{Psi1})$ the previous results are still true if
   $\psi\in
C^1(\overline{\Omega})$, with
 $\psi_{|_{\partial{\Omega}}}\!=\!0$ when
$\emph{K}\!=\!\emph{D}$.
\end{theorem}
\begin{proof}
We will show that under our assumption
$(\ref{condsuaij})-(\ref{bcd})$,
$(\ref{ipotesibijeci})$ on the
coefficients $a_{i,j}$, $b_{i,j}$, $c_j$\, $(i,j=1,2,3)$ of the
linear differential operators $\cal{A},\, \cal{B},\, \cal{C}$
defined in $(\ref{A})$ we can apply the abstract results of Section 4 to
prove locally in time existence and uniqueness of the solution $(u,k)$ to
the identification problem $\PI$.\\
For this purpose let $p\in (3,+\infty)$ and let us
choose the Banach space $X,\,{\wtil X}_1,\,X_1,\,X_2,\, Y,\, Y_1$ appearing in assumptions
$(\textrm{H}1)-(\textrm{H}12)$ according to the rule

\begin{equation}
X=L^p(\Om)\,,\q {\wtil X}_1=W^{1,p}(\Om)\,, \q X_1=W^{1,p}_{\textrm{K}}(\Om)\,,
\q X_2=W_{\textrm{K}}^{2,p}(\Om)\,,
\end{equation}
\begin{equation}
Y=L_{2}^p(0,R)\,, \q Y_1=W_{2}^{1,p}(0,R)\,.
\end{equation}
Since $p\in (3,+\infty)$, reasoning as as in the first part of Section 5 in
\cite{3}, we conclude that $A={\cal A}-\l_0I$ satisfies (H1) -- (H3) in
the sector $\Sigma_{\zeta}$ for some $\l_0 \in \rsp_+$.
\pn
Since assumptions $(\textrm{H}4)-(\textrm{H}6)$ are obviously fulfilled,
we have that $(\textrm{H}1)-(\textrm{H}6)$ hold.
Define now operators $\Phi, \Psi,\, \Psi_1,$
respectively, by $(\ref{Phi1})$, $(\ref{Psi1})$,
$(\ref{psi11ball})$ and operators $E$ and
 $\mathcal{M}$ by
\begin{eqnarray}
\label{EE}Eq(r)=\int_r^{R_2}\!\!q(\xi)d\xi,\qquad\forall\,r\in [0,R],\;
\\[1,8mm]
\label{MM}\qq\mathcal{M}(q,w)(x)=q(|x|)w(x),\qquad\forall\,x\in\Omega,\quad
\end{eqnarray}
Then from H\"older's inequality and the fact that
$p\in (3,+\infty)$ we get\\
\begin{eqnarray}
{\|Eq\|}_{L_{2}^p(0,R)}^p\!\!\!&=&\!\!\!\int_0^R\!r^2{\bigg|
\int_r^Rq(\xi)d\xi\bigg|}^p dr\leqslant
\int_0^R\!r^2{\bigg[\int_0^R\xi^{-2/p}\xi^{2/p}|q(\xi)|
d\xi\bigg]}^p dr\nonumber\\[1,8mm]
\label{EP}&\leqslant&\!\!\!{\|q\|}_{L_{2}^p(0,R)}^p\,\int_0^R\!\!
r^2{\bigg[\int_0^R\!\xi^{^{-{2}/{(p-1)}}}d\xi\bigg]}^{p-1}\!\!\!dr=
\frac{R^{p}}{3}{\Big(\frac{p-1}{p-3}\Big)}^{p-1}
{\|q\|}_{L_{2}^p(0,R)}^p\,.\qq
\end{eqnarray}
Since $D_rEq(r)=-q(r)$ from $(\ref{EP})$ it follows:
\begin{eqnarray}
{\|Eq\|}_{W_{2}^{1,p}(0,R)}\!\!\!&=&\!\!\!
{\Big[{\|Eq\|}_{L_{2}^p(0,R)}^p+
{\|D_rEq\|}_{L_{2}^p(0,R)}^p\Big]}^{1/p}\nonumber\\[1,8mm]
&\leqslant &\!\!\!
{\Big[\frac{R^{^p}}{3}{\Big(\frac{p-1}{p-3}\Big)}^{\!p-1}\!
+1\Big]}^{1/p}{\|q\|}_{L_{2}^p(0,R)}\,.
\end{eqnarray}
Hence $E\!\in\!{\cal{L}}{\big(L_{2}^p(0,R);W_{2}^{1,p}(0,R)\big)}$.
Therefore, by virtue of $(\ref{0.1})$, $(\ref{psinorm})$,
 $(\ref{C1})$ assumption (H7) is
satisfied.\\
Since $p\in (3,+\infty)$ we have the embedding $(\ref{EMB1})$.
Then from the following inequalities,
\begin{eqnarray}
\label{M1}
\|{\cal M}(q,w)\|_{L^p(\Omega)}^p\!\!\! &=&\!\!\!
 \int_\Omega |q(|x|)|^p |w(x)|^p\,dx
\,\le\, \|w\|_{C({\overline{\Omega}})}^p\int_\Omega |q(|x|)|^p\,dx
\nonumber \\[2mm]
&\le&\!\!\! 4\pi \|w\|_{C({\overline{\Omega}})}^p\int_0^R r^{2}|q(r)|^p\,dx
\,\le\, C\|w\|_{W^{1,p}(\Omega)}^p\|q\|_{L_{2}^p(0,R)}^p,\q\;
\end{eqnarray}
we conclude that ${\cal{M}}$ is a bilinear continuous operator from
$L_{2}^p(0,R)\times W^{1,p}(\Om)$ to $L^p(\Om)$. Moreover, using
the embedding $(\ref{EMB2})$ it is an easy task to prove that
 $\cal{M}$ is also continuous from
$W_{2}^{1,p}(0,R)\times L^p(\Om)$ to $L^p(\Om)$ and so (H8) is
satisfied.\\
Then we define $J_1(u_0)$,\,$J_2(u_0)$,\,$J_3(u_0)$ according to formulae
$(\ref{J1}),(\ref{J2}), (\ref{N2})$ and it immediately follows that
assumptions (H9) is satisfied, too.\\
Finally we estimate the vector $(v_0,z_0,z_1,z_2,h_0,q_0)$ in terms of the
 data $(f,u_0,u_1,g_1,g_2)$. Definitions $(\ref{N10})\!-\!(\ref{N0})$ imply
that
\begin{align}
&\hskip 0,5truecm N_1^0(u_1,g_1,f),\;N_3^0(u_0,u_1,g_1,f)\in
C^{\beta}([0,T];L_{2}^{^{p}}(0,R)),&\nonumber\\[1,7mm]
&\hskip 1truecm N_2^0(u_1,g_2,f),\;N_0(u_0,u_1,g_1,g_2,f)\in C^{\beta}([0,T]).&\nonumber
\end{align}
Therefore from $(\ref{h0})$ and $(\ref{q0})$ we deduce
\begin{equation}
\hskip -0,6truecm (h_0,q_0)\in  C^{\beta}([0,T])\times C^{\beta}([0,T];L_{2}^{^{p}}(0,R)),
\end{equation}
whereas from $(\ref{z1z2z3})$, $(\ref{v01})$ and hypotheses
 $(\ref{richiestasuf})\!-\!(\ref{richiestaperA2u0})$ it follows
\begin{align}
&(z_0,z_1,z_2)\in C^{\beta}([0,T];L^{p}(\Omega))\times C^{\beta}([0,T];
W^{1,p}
(\Omega))\times C^{\beta}([0,T];L^{p}(\Omega)),&\\[2,5mm]
&\qq\qq\q v_0\in W_{\textrm{K}}^{2,p}(\Om),\,\q
{\cal{A}}v_0+z_2(0,\cdot)\in W_{\textrm{K}}^{2\beta,p}(\Om)\,.&
\end{align}
Hence assumptions $(\textrm{H}10)\!-\!(\textrm{H}12)$ are also satisfied.
To check condition $(\textrm{H}13)$ first we recall that in this case
the interpolation space ${\cal D}_A(\b,+\infty)$ coincides
with the Besov spaces $B_{\textrm{H,K}}^{2\b,p,\infty}(\Om)\!\equiv
\!{\big(L^p(\Omega), W_{\!\textrm{H,K}}^{2,p}(\Omega)\big)}_{\beta,\infty}$
(cf. $\cite[\textrm{section 4.3.3}]{5}$).
Moreover, we recall that $B_{\textrm{H,K}}^{2\b,p,p}(\Om)=W_{\textrm{H,K}}^{2\b,p}(\Om)$.
Finally, we remind the basic inclusion (cf. \cite[section 4.6.1]{5})
\begin{equation}\label{inclusion}
W^{s,p}(\Om)\hookto B^{s,p,\infty}(\Om)\,,\q\;\textrm{if}\;\,s\notin \mathbb{N}\,.
\end{equation}
Since our function $F$ defined in $(\ref{richiestaperA2u0})$
belongs to $W_{\textrm{H,K}}^{2\b,p}(\Om)$, it is necessarily an
element of $B_{\textrm{H,K}}^{2\b,p,\infty}(\Om)$. Therefore
$(\textrm{H}13)$ is satisfied, too.\ \end{proof}
\textbf{Proof of Theorem \ref{sfera}.} It easily follows from
Theorems \ref{3.1} and \ref{sfera1}.\ $\square$
\begin{remark}\label{su2.36}
\emph{We want here to give some insight into the somewhat involved
condition $(\ref{richiestaperA2u0})$. For this purpose we need to
assume that the functions $a,b,d\in W^{3,\infty}((0,R))$, $c\in
W^{3,\infty}(\Om)$ satisfy the following conditions
\[
b(0)=b'(0)=b''(0)=0,\q d(0)=d'(0)=d''(0)=0,
\]
\[
a'(0)=a''(0)=0,\q D_{x_i}c(0)=D_{x_i}D_{x_j}c(0)=0,\q
i,j=1,\ldots,n.
\]
This implies that the coefficients $a_{i,j}$ belongs
$W^{3,\infty}(\Om)$, $i,j=1,2,3$. Then we observe that function
$k_0$ defined in (3.20) actually belongs to $C^{1+\a}([R_1,R_2])$,
$\a\in (2\b,1)$. It is then an easy task to show the membership of
function $F$ in $W_{\textrm{H,K}}^{2\b,p}(\Om)$, $\b\in (0,1/2)$
under the following regularity assumptions
\begin{itemize}
\item
[{\it{i)\ }}] for any $\rho\in C^{\a}(\ov\Om), \a\in (2\b,1),
w\in W^{2\b,p}(\Om)$, $\rho w \in W^{2\b,p}(\Om)\;$ and satisfies\\
\ the estimate $\|\rho w\|_{W^{2\b,p}(\Om)}\le
C\|\rho\|_{C^{\a}(\ov\Om)} \|w\|_{W^{2\b,p}(\Om)}$\,;
\item
[{\it{ii)\ }}] operator $\Phi$ maps $C^{\a}(\ov\Om)$ into
$C^{\a}([R_1,R_2])$\,.
\end{itemize}
As for as the boundary conditions involved by assumption
$(\textrm{H}13)$ are concerned, we observe that they are missing
when $(\textrm{K})=(\textrm{N})$, while in the remaining case they
are so complicated that we like better not to explicit them and we
limit to list them as
$$F\;\textrm{satisfies boundary conditions (K)}.$$
Of course, when needed, such conditions can be explicitly computed
in terms of the data and function $k_0$ defined in (3.20).}
\end{remark}

\section{The two-dimensional case}
\setcounter{equation}{0}
In this section we deal with
the planar identification problem $\PII$  related to the disk
$\Omega=\{x\in\mathbb{R}^2\!:|x|<R\}$ where $R>0$.\\
Operators $\mathcal{A}$, $\mathcal{B}$, $\mathcal{C}$
are defined by $(\ref{A})$ simply replacing the subscript $3$
with $2$:
\begin{eqnarray}
\label{A1}\mathcal{A}\!=\!\!\sum_{j=1}^{2}D_{x_j}
\big(\sum_{k=1}^{2}a_{j,k}(x)D_{x_k}
\big)\,,\quad\,\mathcal{B}\!=\!\!\sum_{j=1}^{2}D_{x_j}\big(\sum_{k=1}^{2}
b_{j,k}(x)D_{x_k}\big)\,,\quad\,
\mathcal{C}\!=\!\!\sum_{j=1}^{2}c_{j}(x)D_{x_j}\,.
\end{eqnarray}
According to $(\ref{condsuaij})$ for the two-dimensional case, we assume that
the coefficients $a_{i,j}$ of $\cal{A}$ have the following representation
\begin{equation}\label{condsuaij1}
\left\{\begin{array}{lll}
a_{1,1}(x)\!\!\!&=&\!\!\!a(|x|)+\displaystyle\frac{x_2^2[c(x)+b(|x|)]}{|x|^2}
-\displaystyle\frac{x_1^2d(|x|)}{|x|^2},\\[5,0mm]
a_{2,2}(x)\!\!\!&=&\!\!\!a(|x|)+\displaystyle\frac{x_1^2[c(x)+b(|x|)]}{|x|^2}
-\displaystyle\frac{x_2^2d(|x|)}{|x|^2},\\[5,0mm]
a_{1,2}(x)\!\!\!&=&\!\!\! a_{2,1}(x)=-\displaystyle\frac{\,x_1x_2[
b(|x|)+c(x)+d(|x|)]}{|x|^2},
\end{array}\right.
\end{equation}
\pn where the function $a$, $b$, $c$ and $d$ satisfy properties
$(\ref{regular})$, $(\ref{bcd})$. \\
Furthermore we assume that the coefficients of operators
$\mathcal{B},\,\mathcal{C}\,$ satisfy $(\ref{ipotesibijeci})$.\\
In the two-dimensional case, setting $x'=(\cos{\!\varphi},\sin\varphi)$
an example of admissible linear operators $\Phi$
and
$\Psi$ is now the following:
\begin{eqnarray}
\label{Phi12}
\hskip 0,74truecm\Phi [\!\!\!\!\! &v&\!\!\!\!\!](r)\!:=
\int_{\!0}^{2\pi}\!\!\!\!v(rx')d\varphi\,,\qq\\[1,7mm]
\label{Psi12}
\hskip 0,74truecm\Psi[\!\!\!\!\! &v&\!\!\!\!\!]\!:=
\int_{\!0}^{R}\!\! r dr\int_{\!0}^{2\pi}\!\!\!\!\psi(rx')v(rx')\,d\varphi,
\end{eqnarray}
Similarly to $(\ref{tildeA})$,  using $(\ref{condsuaij1})$,
 we  obtain the following
polar representation
for the second order differential operator $\mathcal{A}$:
\begin{eqnarray}
\label{tildeA2}
\widetilde{\mathcal{A}}\!\!\! & = &\!\!\! D_r\big[{h}(r)D_r\big]
\,+\,\frac{{h}(r)D_r}{r}\,+\, \frac{{a}(r)+
{b}(r)}{r^2}D_{\varphi}^2\,+\,\frac{1}{r^2}
D_{\varphi}\big[\,\wtil{c}(r,\varphi)D_{\varphi}\big],
\end{eqnarray}
where $\wtil{c}(r,\varphi)=c(r\cos{\!\varphi},r\sin{\!\varphi})$
and function $h$ is defined in $(\ref{H})$.\\
Working in the Sobolev spaces $W^{k,p}(\Om)$, we will assume
\begin{equation}\label{P2}
p\in (2,+\infty).
\end{equation}
Moreover, our assumptions on operators $\Phi$ and $\Psi$ and the data
will be the same as in $(\ref{primasuPhiePsi})\!-\!
(\ref{richiesteperg2})$ with the spaces $L_2^p(0,R)$ and $W_2^{2,p}(0,R)$
 replaced, respectively, by $L_1^p(0,R)$ and $W_1^{2,p}(0,R)$.
The Banach spaces ${\mathcal{U}}^{\,s,p}(T)$,
${\mathcal{U}}_{\textrm{K}}^{s,p}(T)$ are still
defined by $(\ref{Us})$.
\begin{theorem}\label{sfera2}
Let us suppose that the coefficients $a_{i,j}$ $(i,j=1,2)$ are
represented by $(\ref{condsuaij1})$ and that $(\ref{regular})$,
$(\ref{bcd})$, $(\ref{ipotesibijeci})$,
$(\ref{primasuPhiePsi})\!-\!(\ref{primasuPsi})$,
 $(\ref{P2})$ are fulfilled.
Moreover, assume that the data enjoy the properties $(\ref{richiestasuf})\!
-\!(\ref{richiesteperg2})$ and  satisfy inequalities $(\ref{J0}),
(\ref{J1})$ as well as consistency conditions
$(\ref{DD1})-(\ref{1.19})$, $(\ref{DDV})-(\ref{PSIV1})$.\\
Then there exists $T^{\ast}\in (0,T]$ such that the
identification problem  $\emph\PII\,, \emph{K}\in\{\emph{D,N}\} $,
admits a unique solution $(u,k)\in{\mathcal{U}}^{\,2,p}(
T^{\ast})\times C^{\beta}\big([0,T^{\ast}];W_1^{1,p}(0,R)\big)$
depending
continuously on the data with respect to the norms pointed out in
$(\ref{richiestasuf})\!-\!(\ref{richiesteperg2})$.\\
In the case of the specific operators $\Phi$, $\Psi$ defined as in
$(\ref{Phi12}),\,(\ref{Psi12})$ the previous results are still true if
we assume
  $\psi\in C^1(\overline{\Omega})$ with
${\psi}_{|_{\partial\mbox{}\Om}}\!=\!0$ when $\emph{K}\!=\!\emph{D}$.
\end{theorem}
\begin{lemma}\label{PHIPSI1}
 When $\Phi$ and $\Psi$ are  defined by $(\ref{Phi12})$ and $(\ref{Psi12})$,
respectively, and the coefficients $a_{i,j}$ $(i,j=1,2)$ are
represented by $(\ref{condsuaij1})$, conditions
$(\ref{primasuPhiePsi})\!-\!(\ref{primasuPsi})$ are
 satisfied under assumptions $(\ref{regular})$,
$(\ref{P2})$ and the hypothesis
 $\psi\in C^1(\overline{\Omega})$ with
 ${\psi}_{|_{\partial\mbox{}\Om}}\!=\!0$ when $\emph{K}\!=\!\emph{D}$.
\end{lemma}
\begin{proof}
It is essentially the same as that of Lemma \ref{PHIPSIBAll}.
Therefore we leave it to the reader.\ \end{proof}
For the two-dimensional case the results of Section 5 are still true.
Therefore the proof of Theorem \ref{sfera1} is analogous to the one of Theorem
\ref{sfera}.

\end{document}